\documentclass[11pt,a4paper]{amsart}

\usepackage{amssymb}
\usepackage{amsthm}
\usepackage{amsmath}
\usepackage{amstext}

\numberwithin{equation}{section}

\theoremstyle{plain}
\newtheorem{thm}{Theorem}[section]

\newtheorem{cor}[thm]{Corollary}

\theoremstyle{definition}

\newtheorem*{notation}{Notation}

\theoremstyle{remark}
\newtheorem{rem}[thm]{Remark}

\newcommand{\pg}{\S}

\newcommand{\C}{\mathbb{C}}
\newcommand{\R}{\mathbb{R}}
\newcommand{\K}{\mathbb{K}}

\newcommand{\inv}{^{-1}}
\newcommand{\p}{\partial}
\DeclareMathOperator{\grad}{\p}
\DeclareMathOperator{\gradx}{\p_x}
\newcommand{\infe}{\inf_{\eta \in \K}}
\renewcommand{\v}{\Vec {\mathbf{v}}}
\newcommand{\V}{\Vec {\mathbf{V}}}
\newcommand{\w}{\Vec {\mathbf{w}}}
\newcommand{\of}{\overline f}
\newcommand{\pr}{\mathbb P}

\def\({(\!(}
\def\){)\!)}

\title[$\mu$-constant Deformations of Type $f(x) + t g(x)$]
{Topological Triviality of $\mu$-constant \\ 
Deformations of Type $f(x) + t g(x)$}

\author{Adam Parusi\'nski}

\address{D\'eparement de Math\'ematiques, Universit\'e d'Angers,
   2, bd Lavoisier, 49045 Angers cedex 01, France}

\email{parus@tonton.univ-angers.fr}


\newcommand{\abstracttext}{We show that every $\mu$-constant family of 
isolated hypersurface singularities of type $F(x,t) = f(x) + t g(x)$, where 
$t$ is a parameter, is topologically trivial.  In the proof we construct 
explicitely a vector field trivializing the family. The proof uses only 
the curve selection lemma and hence, for an appropriately translated statement, 
also works over the reals.  Some applications to the study of singularities at 
infinity of complex polynomials are given.}

\begin{document}  
\begin{abstract} \abstracttext \end{abstract} \maketitle

\vspace{ 1.5 truecm}
The main purpose of this note is to show the following result.  

\begin{thm}\label{triv}
Let $\K= \R \text { or } \C$.  
Let $F(x,t)= f(x) +t g(x)$, where $f,g:(\K^n,0) \to (\K,0)$ are germs of 
$\K$-analytic functions and $t\in \K$.  Suppose that in a neighbourhood 
of the origin 
\begin{equation}
|g(x)| \ll  |\grad F(x,t)| \qquad \text {as } \quad (x,t) 
\to g\inv (0) \cap f\inv (0).  
\end{equation}
Then $F$ is topologically trivial along $t$-axis at the origin.  
\end{thm}\medskip

Over the field of complex numbers theorem \ref{triv} can be used to
study some $\mu$-constant families of singularities.  
Let $F(x,t):(\C ^{n+1}\times \C,0)\to (\C,0)$, $F(0,t)\equiv 0$, be a 
$\mu$-constant family of isolated singularities 
$F_t:(\C ^{n+1})\to (\C,0)$.  The answers to the following simple questions 
are still not known 
\begin{enumerate}
\item 
Is any $\mu$-constant deformation topologically trivial?
\item 
Is any $\mu$-constant deformation equimultiple?
\end{enumerate}
Yet, under some additional assumptions the positive answers have been given.  
In particular, the L\^e-Ramanujan theorem \cite{le-ram} 
answers positively the first question 
provided $n\ne 2$.  We note that the proof of this theorem goes beyond
Analytic Geometry and is based on a nontrivial topological fact,
namely the h-cobordism theorem.  

In \cite{trotman} (see also \cite{greuel}, \cite{donald}), using a
simple and elementary argument, 
the author answered positively the second question provided 
$F$ is of the form 
\begin{equation} 
F(x,t) = f(x) + t g(x).
\end{equation}  
We show in corollary \ref{mu-constant} below that the answer to the first
question is also positive if one assumes (0.2).  

Theorem \ref{triv} is proven in section \ref{Triv} below.  The proof,
based on the ideas of Kuo and {\L}ojasiewicz, uses only the curve
selection lemma as a tool.  Therefore, in particular, it would not be  
difficult to extend the results of this paper to any category admitting 
the curve selection lemma such as the subanalatic category or an
o-minimal structure.  

\medskip
\begin{notation}
For a function $F(x,t)$ we denote by $\grad F$ the gradient of $F$ and by 
$\gradx F$ the gradient of $F$ with respect to variables $x$.  Over the field 
of complex numbers $\grad F = 
\overline {(\p F/\p x_1, \ldots,  \p F/\p x_n, \p F/\p t)}$, where 
the bar denotes the complex conjugation.  

$\K$ denotes the field of real or complex numbers. 
\end{notation}
\vskip25pt

\bigskip
\section{{\L}ojasiewicz-type Inequalities}\label{Loj}
\bigskip
 
\begin{thm}\label{loj}
Let $F:(\K^n,0) \to (\K,0)$, $g:(\K^n,0) \to (\K,0)$ be two germs of 
$\K$-analytic functions.  Then there exists a real constant $C$ such that 
for $p\in  F \inv (0)$ and sufficiently close to 
the origin
\begin{equation}
|g(p)| \le C  |p| \infe ( |\eta \grad F(p) + \grad g(p) |) . 
\end{equation}
\end{thm}\medskip

\begin{proof}
Suppose that this is not the case.  Then, by the curve selection lemma, 
there exist a real analytic curve $p(s)$ and a function $\eta (s)$, 
$s\in [0, \varepsilon)$, such that 
\item {(a)} $p(0)=0$;
\item {(b)} $F(p(s)) \equiv 0$, and hence $dF(p(s)) \Dot p(s) \equiv 0$; 
 \begin{equation*}
 |g(p(s))| \gg  |p(s)| |\eta (s) \grad F(p(s)) +  \grad g(p(s)) |, 
\qquad \text { as } \quad s\to 0 . \tag {c}
\end{equation*}  

Since $p(0)=0$ and $g(p(0))=0$ we have asymptotically as $s\to 0$
\begin{equation*}
s  |\Dot p(s)| \sim | p(s)| \qquad \text {and } \quad s 
\frac d {ds} g(p(s)) \sim g(p(s)) . 
\end{equation*}
But
\begin{equation*}
\frac d {ds} g(p(s)) = dg(p(s)) \Dot p(s) = 
(\eta(s) dF + dg) (p(s)) \Dot p(s) . 
\end{equation*}
Hence
\begin{equation*}
|g(p(s))| \sim | s \frac d {ds} g(p(s))| 
\lesssim s |\Dot p(s)| |\eta(s) \grad F(p(s)) + \grad g(p(s)) |, 
\end{equation*}
which contradicts (c).  This ends the proof.    
\end{proof}
\bigskip

\begin{rem}
By standard argument, based again on the curve selection lemma, (1.1)
 implies that in a neighbourhood of the origin 
\begin{equation}
|g(p)| \ll \infe ( |\eta \grad F(p) + \grad g(p) |)   
\qquad \text {as } p \to g\inv (0). 
\end{equation}
Hence, by \cite{loj}, there exists a constant $\alpha>0$ such that 
\begin{equation}
|g(p)| \le |g(p)|^\alpha \infe ( |\eta \grad F(p) + \grad g(p) |)
\end{equation}
which is a variation of the second 
{\L}ojasiewicz Inequality \cite{loj}[$\pg$18].
\end{rem} 
\bigskip

\begin{cor}\label{ag}
Let $F(x,t)= f(x) +t g(x)$, where $f,g:(\K^n,0) \to (\K,0)$ are germs of 
$\K$-analytic functions and $t\in \K$.  Then the following conditions are 
equivalent:  
\medskip
\item  {(i)} \centerline 
{$|g(x)| \ll  |\grad F(x,t)| \qquad \text {as } \quad (x,t) \to (0,0) .$}
\smallskip         
\item {(ii)} \centerline
{$\quad |g(x)| \ll  \inf_\lambda \{|\grad F(x,t)+\lambda \grad g(x) |\}$ } 
as $(x,t)\in F\inv (0)$  and $(x,t) \to (0,0)$.
\end{cor}\medskip

\begin{proof}
(i) $\Longrightarrow$ (ii) follows immediately from theorem \ref{loj}.  
Indeed,  
$\grad F(x,t) = (\gradx f +t \gradx g, g)$, so (i) is equivalent to 
$|g(x)| \ll  |\gradx f +t \gradx g|$ as $(x,t) \to (0,0)$.  Hence 
$|g(x)| \ll   |\grad F(x,t) + \lambda \grad g(x) |$ for $\lambda$ 
small enough, 
say $\lambda < \varepsilon$.  For $\lambda \ge \varepsilon$ the required 
inequality  follows from theorem \ref{loj} by setting $\mu = \lambda\inv$.  

Suppose now that (i) fails.  Then there exists a real analytic curve 
$(x(s), t(s))\to (0,0)$ such that $|g(x(s))| \ge C |\grad F(x(s),t(s))| 
\gg |F(x(s),t(s))|$.  This follows $f(x(s))/g(x(s)) \to 0$ as $s\to 0$.  
Set $\tilde t(s) = - f(x(s))/g(x(s))$.  Then $(x(s), \tilde t(s))\to (0,0)$, 
$F(x(s), \tilde t(s))\equiv 0$ and 
$g(x(s)) \ge C |\grad F(x,t)+\lambda (s) \grad g(x)|$ for 
$\lambda (s) = t(s) - \tilde t(s)$, which contradicts (ii).  
\end{proof} 
\bigskip

\begin{rem}
 (ii) of corollary \ref{ag} is equivalent to the following: suppose that 
$(x,t)\in F\inv (0)$ and $(x,t) \to (0,0)$ then   
\begin{equation*}
\varphi(x) = |g(x)|/ (\inf_\lambda \{|\gradx f(x) + \lambda \gradx g(x)|\}) 
\to 0  .
\end{equation*}
Note that, though $\varphi(x)$ does not depend on $t$, the existence
of $t=t_x$ such that $(x,t)\in F\inv (0)$, $(x,t) \to (0,0)$, 
gives a nontrivial condition on $x$.  
 
If $\gradx f$ and $\gradx g$ are independent then 
\begin{equation}
\inf_\lambda \{|\gradx f(x) + \lambda \gradx g(x)|^2 = 
\frac {|\gradx f|^2 |\gradx g|^2 - |<\gradx f, \gradx g>|^2} {|\gradx g|^2}, 
\end{equation}
and the minimum is attained at 
$\lambda = - \frac {<\gradx f, \gradx g>} {|\gradx g|^2}$.  
Hence (ii) of corollary \ref{ag} is equivalent to  
\begin{equation} 
\frac {|g(x)|^2 |\gradx g|^2} {\Delta }  \to 0  ,
\end{equation}
where $\Delta = |\gradx f|^2 |\gradx g|^2 - |<\gradx g,\gradx f>|^2$.  
\end{rem} 
\bigskip

\bigskip
\section{Topological Triviality}\label{Triv}\bigskip

\begin{proof} [Proof of theorem \ref{triv}]
First we show how to trivialize the zero set $X = F\inv (0)$.  For this 
we construct a vector field $\v$ on $X$ whose $t$-coordinate is $\p_t$ and
which restricted to the t-axis is precisely $\p_t$.   
 Moreover, for the reason we explain below,  we require 
that $\v$ is tangent to the levels of $g$.  
If one looks for such a vector field of the form 
\begin{equation*}
\v =   \p_t + a \gradx f + b \gradx g 
\end{equation*}
the condition of tangency to the levels of $F$ and $g$ impose that 
on $X'= X\setminus g\inv (0)$ 
\begin{equation}
\v =   \p_t - {\frac {\bar g |\gradx g|^2} \Delta}  \gradx f + 
{\frac {\bar g <\gradx f,\gradx g>} \Delta} \gradx g ,  
\end{equation}
where again $\Delta = |\gradx f|^2 |\gradx g|^2 - |<\gradx g,\gradx f>|^2$.  

We claim that 
\begin{equation}
|\v - \p_t| \to 0  \qquad \text {as } \quad (x,t) \to g\inv (0) , F(x,t)=0 . 
\end{equation}
It suffices to show it on any real analytic curve  
$(x(s),t(s))\to (x_0,t_0) \in g\inv (0)\cap F\inv(0)$.  We may assume 
$(x_0,t_0) =(0,0)$.  Since 
$|\v - \p_t|^2 = \frac {|g(x)|^2 |\gradx g|^2} {\Delta }$, (2.2) follows 
from Corollary \ref{ag} and (1.5).  

To show that the flow of $\v$ is continuous we use $g$ as a control function.  
First note that, by (1.4), $X' =X \setminus g\inv (0)$ is nonsingular.  
By (1.8),  
$\{\Delta =0\} \subset \{g=0\}$, and hence $\v$ is smooth on $X'$.  Moreover 
the integral curves of $\v|_{X'}$ cannot fall on $Y=X\cap g\inv (0)$ since 
$\v$ is tangent to the levels of $g$.  Extend $\v$ on $Y$ by setting 
$\v|_Y \equiv \p_t$.  Let $\Phi(p,s)$ be the flow of $\v$ and let 
$p_0 = (x_0,t_0)\in Y$, $p_1 = (x_1,t_1) \in X'$.  By (2.2) 
\begin{equation}
|\v(\Phi (p_0,s)) - \v(\Phi (p_1,s))| \le \varepsilon (g(\Phi (p_1,s))) = 
\varepsilon  (g(p_1)) , 
\end{equation}
and $\varepsilon (\eta)\to 0$ as $\eta \to 0$.  Hence 
\begin{equation}
|(\Phi (p_0,s)) - (\Phi (p_1,s))| \le |p_0-p_1| + 
|s| \varepsilon (g(p_1)) .
\end{equation}
This shows the continuity of $\Phi$ and ends the first part of the proof 
of theorem, the triviality of $X$ along the $t$-axis.  

To trivialize $F$ as a function we use the Kuo vector field \cite{kuo} 
\begin{equation}
\w =   \p_t - {\frac {\bar g} {|\gradx F|^2} }  \gradx F  .  
\end{equation}
$\w$ is tangent to the levels of $F$ and hence, by (0.1), 
one may show that the flow of $\w$ is continuous in the 
complement of $X=F\inv (0)$.   Unfortunately this argument does not work 
on $X$ and the integral curves of $\w|_X$ may disappear at 
the singular locus of $F$.  To overcome this difficulty we "glue" $\v$ and 
$\w$.  Fix a sufficiently small neighbourhood $\mathcal U$ af the origin in 
$\K^n \times \K$.  Let $\mathcal V_1$ be a neighbourhood of $X'=X\setminus Y$ 
in $\mathcal U \setminus Y$ such that (1.5) still holds on $\mathcal V_1$.  
Fix a partition of unity $\rho_1, \rho_2$ subordonated to the covering 
$\mathcal V_1,\mathcal V _2 = \mathcal U\setminus X$ of 
$\mathcal U \setminus Y$ and put
\begin{equation*}
\V =   \rho_1 \v + \rho_2 \w  .  
\end{equation*}

Clearly $\V$ is tangent to the levels of $F$ on $\mathcal U \setminus Y$ 
and to the levels of $g$ on $X \setminus Y$.  Hence its flow is continuous on 
$\mathcal U \setminus Y$.  Moreover 
\begin{equation*}
|\V - \p_t| \le  \rho_1 |\v- \p_t| + \rho_2 |\w- \p_t|  ,  
\end{equation*}
so it goes to $0$ as $(x,t)\to Y$.  Now the continuity of the flow of 
$\V$ follows by the same argument as the continuity of the flow of $\v$.

This ends the proof of theorem.  
\end{proof}

\medskip
\begin{cor}\label{mu-constant} 
Let $f,g:(\C^n,0) \to (\C,0)$ be such that $F(x,t)= f(x) +t g(x)$ is
 a $\mu$-constant family of isolated singularites.  
Then $F$ is topologically trivial.   
\end{cor}\smallskip

\begin{proof}
Let $L$ denote the t-axis.  
By \cite{le-saito}  $\mu$-constancy is equivalent to 
\begin{equation}
|\p_t F| \ll  |\grad F(x,t)| ,  
\end{equation}
as $(x,t) \to L $.
Hence, by theorem \ref {triv} it suffices to show (2.6) 
for $(x,t)\to (x_0,t_0)$,  where $g(x_0)=0$ and $x_0\ne 0$.  
But this is obvious since, by our hypothesis, $\gradx F(x_0,t_0) \ne 0$.  
\end{proof}
\bigskip

\begin{rem}
Suppose that $F(x,t)= f(x) +t g(x)$ is a family of isolated singularities 
and that in a neighbourhood of the origin 
\begin{equation}
|g(x)| \le  C |x| |\grad F(x,t)| .  
\end{equation}
Then, as shown in  \cite{kuo}, the vector field $\w$ of (2.5) trivializes 
$F$ along t-axis.  The flow of $\w$ is continuous since 
$\rho (x,t) = |x|$ changes slowly along the trajectories of $\w$.  That is,  
for a trajectory $(x(s), s)$
\begin{equation}
|\frac {d\rho }{ds}| \le  C \rho ,  
\end{equation}
which implies 
\begin{equation}
\rho (x(0),0) e^{-C|s|} \le |\rho (x(s),s)| \le \rho (x(0),0) e^{C|s|},  
\end{equation}
see \cite{kuo} for the details.  

It is not difficult to see that 
(2.8) and (2.9) are also satisfied on 
the trajectories of vector field $\v$ of (2.1), provided (2.7) holds.  
  Indeed, the proof of theorem \ref{loj} actually gives
\begin{equation}
|g(p)| \le C \infe (|x| |\eta \gradx F(p) +\grad g(p)| + 
|t|| \eta \partial F/\partial t (p)|) . 
\end{equation}

This allows one to show that (2.7) is equivalent to 
\begin{equation*}
|g(x)| \le C |x| \inf_\lambda \{|\grad F(x,t)+\lambda \grad g(x) |\}  
\quad \text { on } F\inv (0).
\end{equation*}
Hence $|\v(x,t) - \p_t| \le C |x|$ which easily implies (2.8) and (2.9). 
\end{rem}
\vskip20pt

\bigskip
\section{Singularities at infinity of complex polynomials}\label{infty}
\bigskip

Let $f(x_1,\ldots ,x_n)$ be a complex polynomial of degree $d$.     
Let $\tilde f(x_0, x_1, \ldots, x_n)$ denote the homogenization of $f$.  
Set 
\begin{equation*}
X = \{ (x,t)\in \pr ^n \times \C | 
\, F(x,t) = \widetilde f(x) - tx_0^d  =0\}
\end{equation*}   
and let $\of:X\to \C$ be induced by the projection on the second factor.  
Then $\of$ is the family of projective closures of fibres of $f$.  

Let $H_{\infty}= \{x_0=0\} \subset \pr ^n$ be the hyperplane at infinity and 
let $X_{\infty}=X\cap (H_\infty \times \C)$.  
$\of$ can be used to trivialize a non-proper function $f$ as follows. 
Fix $t_0\in \C$.  If $\of$ is topologically trivial over a neighbourhood of 
$t_0$ and the trivialization preserves $X_{\infty}$, then clearly 
$f$ is topologically trivial over the same neighbourhood of $t_0$.  
Thus one may use Stratification Theory in order to trivialize $f$.  

On the other hand, by a simple argument, see for instance \cite{parus1}, 
the following condition on the asymptotic behaviour of $\grad f$ at 
infinity assures the topological triviality of $f$ over a neighbourhood of 
$t_0$: $\grad f$ does not vanish on $f\inv (t_0)$ and 
\begin{equation}
|x| |\grad f(x)| \ge \delta > 0  ,
\end{equation}  
for all $x$ such that $|x|\to \infty$ and $f(x) \to t_0$.  
Following \cite{pham} the condition (3.1) is called Malgrange's Condition.  

These two approaches are related by the following result.  

\bigskip
\begin{thm}\label{banach} \cite{parus2, tibar}  \par 
Let $p_0\in X_\infty$, $t_0=f(p_0)$.   Then the following conditions are 
equivalent: 
\item {(i)}
      $\of$ has no vanishing cycles in a neighbourhood of $p_0$. 
\item {(ii)} 
      (3.1) holds as $(x,f(x))\to p_0$ in $\pr^n\times \C$.  
\end{thm}

We show below that the proofs of theorem \ref{banach} given in 
\cite{parus2, tibar} can be considerably simplified using a version 
of corollary \ref{ag}.  
For this we first sketch the main steps of the proof of \cite{parus2}.  

We work locally at $p_0\in X_\infty$.   
We may assume that $p_0 = ((0:0:\ldots :0:1),0) \in \pr ^n\times \C $, 
so that $y_0 = x_n^{-1}$, $y_i = x_i/ x_n$ for $i=1, \ldots, n-1$, and $t$, 
 form a local system of coordinates at $p_0$.  In this new coordinate
 system  $X$ is defined by
\begin{equation*}
F(y_0, y_1, \ldots, y_{n-1},t) = 
\tilde f(y_0, y_1, \ldots, y_{n-1},1)  - t y_0^d = 0 .
\end{equation*} 

It is easy to see that (ii) of theorem \ref{banach} is equivalent to 
\begin{equation}
|\partial F/\partial t (p)| \le C |\partial F/\partial y_1, 
\ldots, \partial F/\partial y_{n-1})(p)| \, ,
\end{equation} 
for $p\in X$ and close to $p_0$.  On the other hand, less trivial but  
standard arguments, see e.g. \cite{parus2}[Definition-Proposition 1], show 
that (i) of theorem \ref{banach} is equivalent to the fact that $(p_0,dt)$ 
does not belong to the characteristic variety of $X$, which is equivalent to 
 \begin{equation}
|\partial F/\partial t (p)| \le C | 
 \partial F/\partial y_0, 
\ldots, \partial F/\partial y_{n-1})(p)| \, ,
\end{equation} 
in a neighbourhood of $p_0$.   

The proofs of  (3.2) $\Longleftrightarrow$ (3.3)  given in \cite{parus2,tibar} 
are based on the results of \cite{BMM}.  But it may be easily shown using 
theorem \ref{loj}.  Indeed, take $g(y) = y_0^d$.  Then 
 \begin{equation*}
  \inf_\lambda \{|\grad F(y,t)+\lambda \grad g(y) |\} =  
|\partial F/\partial y_1, \ldots, \partial F/\partial y_{n-1})(p)| ,
\end{equation*} 
and (3.2) $\Longleftrightarrow$ (3.3) is a variation of corollary \ref{ag} 
with the same proof.  

It seems to be particularly surprising that unlike the previous proofs 
the one presented above does 
not use the complex structure at all.


\end{document}